\documentclass[12pt]{amsart}
\usepackage{xypic}
\xyoption{all}
\usepackage{amscd,amsmath,latexsym,amssymb}
\usepackage[mathcal]{euscript}

\usepackage[all]{xy}
\CompileMatrices
\numberwithin{equation}{section}

\newtheorem{thm}{Theorem}[section]
\newtheorem{prop}[thm]{Proposition}
\newtheorem{lem}[thm]{Lemma}
\newtheorem{cor}[thm]{Corollary}

\theoremstyle{definition}
\newtheorem{defn}[thm]{Definition}
\newtheorem{remark}[thm]{Remark}
\newtheorem{example}[thm]{Example}

\newcommand{\F}{\mathbb F}

\newcommand{\Sym}{\operatorname{S}}
\newcommand{\Char}{\operatorname{char}}
\newcommand{\Soc}{\operatorname{Socle}}

\begin{document}
\title[Bundles, Cohomology and Truncated Symmetric Polynomials]{Bundles, 
Cohomology and Truncated Symmetric Polynomials}

\author[Alejandro Adem]{Alejandro Adem$^{*}$}
\address{Department of Mathematics
University of British Columbia, Vancouver BC V6T 1Z2, Canada}
\email{adem@math.ubc.ca}
\thanks{$^{*}$Partially supported by NSERC and NSF} 
% The first author is grateful
% to W. Dwyer for his helpful comments.}

\author[Zinovy Reichstein]{Zinovy Reichstein$^{**}$}
\address{Department of Mathematics
University of British Columbia, Vancouver BC V6T 1Z2, Canada}
\email{reichst@math.ubc.ca}
\thanks{$^{**}$Partially supported by NSERC}

\date{\today}

\begin{abstract}
The cohomology of the classifying space $BU(n)$ of the unitary groups
can be identified with the
the ring of symmetric polynomials on $n$ variables by restricting
to the cohomology of $BT$, where $T\subset U(n)$ is a maximal torus.
In this paper we explore the situation where
$BT = (\mathbb CP^{\infty})^n$ is replaced by a
product of finite dimensional
projective spaces $(\mathbb CP^d)^n$, fitting into
an associated bundle
$$U(n)\times_T (\mathbb S^{2d+1})^n\to (\mathbb CP^d)^n\to BU(n).$$
We establish a purely algebraic version of this problem by exhibiting 
an explicit system of generators for the ideal of truncated 
symmetric polynomials. We use this algebraic result to give 
a precise descriptions of the kernel of the homomorphism 
in cohomology induced by the natural map $(\mathbb CP^d)^n\to BU(n)$. 
We also calculate the
cohomology of the
homotopy fiber of the natural map $E \Sym_n\times_{\Sym_n}(\mathbb
CP^d)^n\to BU(n)$.
\end{abstract}


% 05E05 (1991-now)  Symmetric functions
% 20J06 (1980-now)  Cohomology of groups
% 57S25 (1980-now)  Groups acting on specific manifolds 
% 55R35 (1980-now)  Classifying spaces of groups and ${H}$-spaces

\maketitle

\tableofcontents

\section{Introduction}
One of the nicest calculations in algebraic topology is that of the
cohomology of the classifying space $BU(n)$ of the unitary groups as
the ring of symmetric polynomials on $n$ variables (see \cite{Bor}). 
In
fact the restriction map identifies $H^*(BU(n),\mathbb Z)$ with the
invariants in the cohomology of the classifying space $BT$ of a
maximal torus under the action of the Weyl group $\Sym_n$. This leads
to a beautiful description of the cohomology of the flag manifold
$U(n)/T$ and more specifically a detailed understanding of the
fibration $U(n)/T\to BT\to BU(n)$.

In this paper we explore the situation where 
$BT = (\mathbb CP^{\infty})^n$ is replaced by a
product of finite dimensional
projective spaces $(\mathbb CP^d)^n$, fitting into
an associated bundle
$$U(n)\times_T (\mathbb S^{2d+1})^n\to (\mathbb CP^d)^n\to BU(n).$$
This requires an analysis of truncated symmetric invariants and in
particular a precise description of the kernel $I(n,d)$ of the
algebra surjection $H^*(BU(n),\mathbb F)\to H^*((\mathbb CP^d)^n,
\mathbb F)^{\Sym_n}$.
The purely algebraic version of this 
problem is 
studied in~\S\ref{sect5} and ~\S\ref{sect.reg-seq}.
In particular, Theorem~\ref{thm1}
allows us to exhibit 
an explicit set of generators for $I(n, d)$ as follows.

% \newpage

\begin{thm} \label{thm1.1} Let $\F$ be a field and $I(n,d)$ be 
the kernel of the map $H^*(BU(n),\F)\to H^*((\mathbb CP^d)^n,\F)$. 

\smallskip
(a) If $n!$ is invertible in $\mathbb F$ then $I(n, d)$
is generated by the elements $P_{d+1}, P_{d+2}, \dots, P_{d+n}$

\smallskip
(b) If $n < 2\Char(\F) - 1$ then $I(n, d)$ is 
generated by $P_{d+1}, P_{d+2}, \dots, P_{d+n}$ and
$P_{\underbrace{\text{\tiny $d+1, \ldots, d+1$}}_{\text{$p$ times}}}$.
\end{thm}

%\begin{thm} \label{thm1.1}
%Suppose $\F$ is a field of characteristic $p > 0$ and
%$k$ is an integer such that $p^k\le n < p^{k+1}$. Then
%the kernel $I(n,d)$ of the map $H^*(BU(n),\F)\to H^*((\mathbb
%CP^d)^n,\F)$ is generated by the following $n+k$ elements:
%\begin{itemize}
%\item $P_{d+i}$, where $1\le i\le n$, and % $|P_j| = 2j$,
%\item $P_{\underbrace{\text{\tiny $d+1, \ldots, d+1$}}_{\text{\tiny $p^i$ times}}}$,
%where $1\le i\le k$. 
%\end{itemize}
%\end{thm}
\noindent For the definition of $P_{d + i}$ and $P_{\underbrace{\text{\tiny $d+1, \ldots,
d+1$}}_{\text{\tiny $p^i$ times}}}$, see~\S\ref{sect5}. 
Note that the degree of $P_{d + i}$ is $2(d + i)$ and the degree of  
$P_{\underbrace{\text{\tiny $d+1, \ldots,
d+1$}}_{\text{\tiny $p$ times}}}$ is $2p(d+1)$.

If $n!$ is invertible in a field $\F$, then we show that
the elements $P_{d+i}$, $1\le i\le n$ form a 
generating regular sequence for $I(n, d)$. In contrast,
using Theorem~\ref{thm2} we show that in most other cases $I(n,d)$
cannot be generated by a regular sequence: 

\begin{thm}
If $n\ge \Char(\F)>0$ and $d>1$, then $I(n,d)$ cannot be generated
by a regular sequence. 
\end{thm}

There is a free action of $\Sym_n$ on the fiber space
$W(n,d)=U(n)\times_T (\mathbb S^{2d+1})^n$ % \to (\mathbb CP^d)^n$
arising from the normalizer of the maximal torus in $U(n)$. The
orbit space $X(n,d)$ can be realized as the fiber of the natural map
$E\Sym_n\times_{\Sym_n}(\mathbb CP^d)^n\to BU(n)$. Our algebraic
calculations allow us to calculate the cohomology of this space in
good characteristic.

\begin{thm}
If $\F$ is a field where $n!$ is invertible, then the cohomology of
$X(n,d)$ is an exterior algebra on $n$ generators
$$H^*(X(n,d), \F)\cong
\Lambda_\F (E_{d+1}, \dots, E_{d+n})$$ 
where $E_j$ is a cohomology
class in dimension $2j-1$.
\end{thm}

This has an interesting computational consequence.

\begin{thm}
For any field $\F$ of coefficients, the Serre spectral sequence for
the fibration $(\mathbb S^{2d+1})^n\to W(n,d)\to U(n)/T$ collapses
at $E_2$ if and only if $d\ge n-1$. Consequently,
we obtain an additive calculation
$$H^*(W(n,d),\F)\cong H^*(U(n)/T)\otimes H^*((\mathbb S^{2d+1})^n,\F)$$
whenever $d \ge n - 1$.
In particular if $n!$ is invertible in $\F$, then
$$H^*(X(n,d),\F)\cong
[H^*(U(n)/T)\otimes H^*((\mathbb S^{2d+1})^n,\F)]^{\Sym_n} \cong
\Lambda_\F (E_{d+1}, \dots, E_{d+n}) \, .$$
\end{thm}

These results follow from a general theorem about the cohomology of
fibrations which, although ``classical'' in nature, seems to be new.

\begin{thm}
Let $\F$ be a field and
let $\pi: E\to B$ denote a fibration with fiber $F$ of finite type
such that $B$ is simply connected. Assume 
\begin{itemize}
\item $H^*(B,\F)$ is a polynomial algebra on $n$ even
dimensional generators,
\item $\pi^*: H^*(B,\F)\to H^*(E,\F)$ is surjective,
\item the kernel of $\pi^*$ is generated by a regular
sequence $u_1,\dots, u_n$, where $|u_i|= 2r_i$. 
% (cf. Definition~\ref{def.regular-sequence}),
\end{itemize}
Then $H^*(F,\F)$ is an exterior algebra on $n$ odd dimensional
generators $e_1, \dots , e_n$, where $|e_i| = 2r_i -1$.
\end{thm}

%%%%%%%%%%%%%%%%%%
%We also provide information about the analogous situation for the
%orthogonal groups, namely about the map induced in mod $2$
%cohomology by the
%group inclusion of diagonal matrices $(\mathbb Z/2)^n\to O(n)$ with
%Weyl group $\Sym_n$.
%%%%%%%%%%%%%%%%%%%
It is natural to ask whether the results of this paper 
can be extended to compact Lie groups, other than $U(n)$. 
We thus conclude this introduction with the following open problem.

\medskip

\noindent\textbf{Problem}: 
Let $G$ be a compact Lie group with maximal torus $T$ 
of rank $n$ and Weyl group $W$. Describe generators for the kernel
$I_G(n,d)$ of the natural map $H^*(BG,\F) \to
H^*((\mathbb CP^d)^n,\F)$ 
and use this to describe the cohomology
of the homotopy fiber of $(\mathbb CP^d)^n\to BG$ when $|W|$ is
invertible in $\mathbb F$.

\medskip
Theorems~\ref{thm1}(a) and~\ref{thm2}(a) have been independently proved 
in a recent preprint~\cite{CKW} by A. Conca, Ch. Krattenthaler, 
J. Watanabe. We are grateful to J. Weyman for bringing 
this preprint to our attention.

\section{Bundles and symmetric invariants}

A classical computation in algebraic topology is that of the
cohomology of the classifying space $BU(n)$ where $U(n)$ is
the unitary group of $n\times n$ matrices.
We briefly recall how that goes; details can be 
found, e.g., in the survey paper \cite{Bor} by A. Borel. 
 Let $T=(\mathbb S^1)^n\subset U(n)$ denote the
maximal torus of diagonal matrices
in $U(n)$; its classifying space is
$BT=(\mathbb CP^\infty)^n$. The inclusion $T\subset U(n)$
induces a map between the cohomology of $BU(n)$ and the cohomology
of $BT$. Note that the normalizer $NT$ of the torus is a
wreath product $\mathbb S^1 \wr \Sym_n$, where the symmetric group
$\Sym_n$ acts by permuting the $n$ diagonal entries. Thus the
Weyl group $NT/T$
is the symmetric group $\Sym_n$. Recall that
$H^*(BT,\mathbb Z)\cong \mathbb Z [x_1,\dots , x_n]$ where
the $x_1,\dots ,x_n$ are 2-dimensional generators.

\begin {thm}
The inclusion $T\subset U(n)$ induces an inclusion in cohomology with
image the ring of symmetric invariants in the graded polynomial
algebra

$$H^*(BU(n),\mathbb Z) \cong H^*(BT,\mathbb Z)^{\Sym_n}
= \mathbb Z [x_1, \dots , x_n]^{\Sym_n}$$
where the action of $\Sym_n$ arises from that of the Weyl group.
\qed
\end{thm}

Now recall that the complex projective space $\mathbb CP^d$
is a natural subspace of $\mathbb CP^\infty$; this induces a map
$$\tilde{F}(n,d): (\mathbb CP^d)^n \to BT\to BU(n).$$
The permutation matrices $\Sym_n\subset U(n)$ act via
conjugation on $U(n)$; this restricts to an action on the
diagonal maximal torus $T$ which permutes the factors.
Applying the classifying space functor yields actions of
$\Sym_n$ on $BT$ and $BU(n)$ which make the map
$\tilde{F}(n,d)$ equivariant.
Note however that the conjugation
action on $U(n)$ is homotopic to the identity on $BU(n)$.
We conclude that $\tilde{F}(n,d)$ induces the 
natural map 
$$\tilde{F}(n,d)^*: H^*(BU(n),\mathbb Z)\cong
\mathbb Z [x_1, \dots , x_n]^{\Sym_n}
\to \mathbb Z [x_1, \dots , x_n]/(x_1^{d+1}, \dots
x_n^{d+1})$$
in integral cohomology whose
image is precisely the ring of truncated symmetric invariants.
We should also note that the map $\tilde{F}(n,d)$ is
(up to homotopy) the classifying map for the $n$--fold
product of the canonical complex line bundle over
$\mathbb CP^d$.

To make this effective geometrically, we need to describe the
map $\tilde{F}(n,d)$ explicitly as a fibration. The space
$(\mathbb CP^d)^n$ is a quotient of $(\mathbb S^{2d+1})^n$
by the free action of the maximal torus $T$. Using a standard
induction construction we can view our map as a fibration
which lies over the classical
fibration connecting $U(n)/T$, $BT$ and $BU(n)$.
Indeed, the following commutative diagram has fibrations
in its rows and columns:
\[  \xymatrix{
& (\mathbb S^{2d+1})^n\ar@{=}[r] \ar[d] & (\mathbb S^{2d+1})^n \ar[d] \\
W(n,d)\ar@{=}[r] & U(n)\times_T(\mathbb S^{2d+1})^n \ar[r] \ar[d] &
(\mathbb CP^d)^n \ar[r]^{\tilde{F}(n,d)} \ar[d] & BU(n)\ar@{=}[d] \\
& U(n)/T\ar[r]  & BT \ar[r]
& BU(n)\\
}
\]
Note that we also have a bundle
\[
U(n)\to U(n)\times_T (\mathbb S^{2d+1})^n\to
(\mathbb CP^d)^n
\]
and its classifying map is $\tilde{F}(n,d)$.

In some of our applications it will also
make
sense to take a quotient by the action of
the symmetric group $\Sym_n$. However for technical
reasons this requires
taking a \textsl{homotopy orbit space} which
we now define.

\begin{defn}
Let $G$ denote a compact Lie group acting on a
space $X$, its homotopy orbit space $X_{hG}$
is defined
as the quotient of the product space
$EG\times X$ by the diagonal
$G$--action, where $EG$ is the universal $G$--space.
\end{defn}

\begin{remark}
It should be noted that if $G$ is a finite group,
$X$ is a $G$--space and $|G|$ is invertible in 
the coefficients, then the natural projection
$X_{hG}\to X/G$ induces an isomorphism in cohomology
(this follows from the Vietoris-Begle theorem).
Hence for example if $|G|$ is invertible in a coefficient
field $\F$, then
$H^*(X_{hG}, \F)\cong H^*(X,\F)^G$ (the algebra of invariants).
\end{remark}

In our context, the symmetric group $\Sym_n$ acts by permuting
the factors in $(\mathbb CP^d)^n$ and we can consider
the associated homotopy orbit space
$$(\mathbb CP^d)^n_{h\Sym_n}=
E\Sym_n\times_{\Sym_n} (\mathbb CP^d)^n.$$
More precisely,
the map $BT\to BU(n)$ naturally factors through the
classifying space of the normalizer $NT$, as we have
$T\subset NT\subset U(n)$. The space $BNT$ can be
identified with $BT_{h\Sym_n}=(\mathbb CP^\infty)^n_{h\Sym_n}$,
where $\Sym_n$ acts by permuting the factors, as before. This homotopy
orbit space restricts to the truncated projective spaces,
yielding a map
\[ F(n,d): (\mathbb CP^d)^n_{h\Sym_n} \to BU(n) \, , \]
which is surjective in rational cohomology. We would also like
to describe this map as a fibration.
\medskip

The map
$(\mathbb CP^d)^n \to BT$ is an $\Sym_n$-equivariant
fibration, with fiber $(\mathbb S^{2d+1})^n$. This arises
from the free $T$--action on the product of spheres,
which extends in the usual way to an action of the
semidirect product $NT$.
If we take
homotopy orbit spaces we obtain a fibration sequence
\[
(\mathbb S^{2d+1})^n \to (\mathbb S^{2d+1})^n_{hNT}\to
BNT \, . 
\]
Dividing out by the free $T$--action we can identify
$(\mathbb S^{2d+1})^n_{hNT}\simeq  (\mathbb CP^d)^n_{h\Sym_n}$.
This makes the fiber
of the map $(\mathbb CP^d)^n_{h\Sym_n}\to BNT$ very explicit.
As before, in order to describe the fibration
with
target $BU(n)$, it suffices to
induce
up the action on the fiber
to a $U(n)$--action by taking the balanced product
$Z=U(n)\times_{NT}(\mathbb S^{2d+1})^n$.
This yields a fibration sequence
\[
Z\to Z_{hU(n)}
\to BU(n) \, .
\]
Note that
$$Z_{hU(n)}\simeq EU(n)\times_{NT}(\mathbb S^{2d+1})^n
\simeq (\mathbb S^{2d+1})^n_{hNT}\simeq
(\mathbb CP^d)^n_{h\Sym_n}$$
where the last equivalence follows from taking quotients
by the free $T$--action, as before. Our discussion 
is summarized in the following
diagram of fibrations, analogous to the non--equivariant
situation:

\[
\xymatrix{
& (\mathbb S^{2d+1})^n\ar@{=}[r] \ar[d] & (\mathbb S^{2d+1})^n \ar[d] \\
X(n,d)\ar@{=}[r] & U(n)\times_{NT}(\mathbb S^{2d+1})^n \ar[r] \ar[d] &
E\Sym_n\times_{\Sym_n}(\mathbb CP^d)^n~~ \ar[r]^{F(n,d)} \ar[d] & BU(n)\ar@{=}[d] \\
& U(n)/NT\ar[r]  & BNT \ar[r]
& BU(n)\\
}
\]

Hence we have

\begin{prop}
Up to homotopy the map $\tilde{F}(n,d):(\mathbb CP^d)^n\to BU(n)$
is a fibration with fiber the compact simply connected
manifold
$$W(n,d)=U(n)\times_T(\mathbb S^{2d+1})^n$$
of dimension equal to
$n(n +2d)$. There is a free $\Sym_n$--action on this
manifold, and its quotient
$$X(n,d)= U(n)\times_{NT}(\mathbb S^{2d+1})^n$$
is homotopy equivalent to the fiber
of $F(n,d):(\mathbb CP^d)^n_{h\Sym_n}\to BU(n)$.
\qed
\end{prop}
\begin{remark}
Note that there are fibrations
\[
(\mathbb S^{2d+1})^n \to X(n,d) 
\to U(n)/NT
\]
and
\[
U(n) \to X(n,d)
\to (\mathbb CP^d)^n_{h\Sym_n}
\]
where the second one is
obtained from pulling back the universal $U(n)$ bundle
over $BU(n)$ using $F(n,d)$.
\end{remark}

One of our main results in this paper will be to calculate the
cohomology of the fibers $W(n,d)$ and $X(n,d)$ associated to the
fibrations $\tilde{F}(n,d)$ and $F(n,d)$ respectively.

\section{Cohomology calculations when $n!$ is invertible}
Our standing assumption in this section (unless stated otherwise)
will be that $\F$ is a field such that $n!$ is invertible in $\F$, 
and cohomology will be computed with $\F$--coefficients. A good example
is the field $\mathbb Q$ of rational numbers. In this situation we have
$H^*(X(n,d),\F) \cong H^*(W(n,d),\F)^{\Sym_n}$; it is this cohomology
algebra that we will be most interested in.

We begin by considering the limit case $d=\infty$. In this
case $X(n,\infty)=U(n)/NT$ and we are looking at the
classical fibration
\[
 U(n)/NT \to BNT
\to BU(n)
\]

\begin{prop}
The map $BNT\to BU(n)$ induces an isomorphism in
cohomology and $U(n)/NT$ is $\F$--acyclic.
\end{prop}
\begin{proof}
Indeed, both maps in the sequence
$$H^*(BU(n),\F)\to H^*(BNT,\F)\to H^*(BT,\F)^{\Sym_n}$$
are isomorphisms. Since $BU(n)$ is simply connected,
this can only happen if $U(n)/NT$ is acyclic.
\end{proof}
Note that this computation is very different from what the
cohomology of
the flag manifold $U(n)/T$ looks like; when we divide out by
the action of the symmetric group all the reduced cohomology
vanishes.

We now consider the unstable case of this result, namely
when $d$ is finite. This is considerably more interesting,
as we know that the cohomology must be non--trivial.
This calculation will be a special case of a more
general result about the cohomology of fibrations.

\begin{thm}\label{EMSS}
Let $\pi: E\to B$ denote a fibration with fiber $F$ of
finite type such
that $B$ is simply connected and
\begin{itemize}
\item $H^*(B,\F)$ is a polynomial algebra on $n$ even
dimensional generators,
\item $\pi^*: H^*(B,\F)\to H^*(E,\F)$ is surjective,
\item the kernel of $\pi^*$ is generated by a regular
sequence $u_1,\dots, u_n$ where $|u_i|= 2r_i$.
\end{itemize}
Then $H^*(F,\F)$ is an exterior algebra on $n$ odd
dimensional generators $e_1, \dots , e_n$, where
$|e_i| = 2r_i -1$.
\end{thm}

\begin{proof}
The cohomology of the fiber $F$ in a fibration
\[
F\to E
\to B
\]
can be studied using the Eilenberg--Moore spectral sequence.
We refer the reader to \cite{mccleary}, Chapter VIII
for details.
It has the form:
$$E_2^{*,*} =
Tor_{H^*(B,\F)}(\F, H^*(E,\F)).$$
On the other hand, the hypotheses imply that
$$H^*(E,\F)\cong H^*(B,\F)/(u_1, \dots u_n)$$
where $u_1, \dots , u_n$ form a regular sequence
of maximal length
in $H^*(B,\F)$, a polynomial algebra on $n$
even dimensional generators. In other words
the cohomology of $B$ is free and finitely
generated over $\F [u_1,\dots , u_n]$.
%$$(\F[x_1,\dots , x_n]/(x_1^{d+1}, \dots ,x_n^{d+1}))^{\Sym_n}
%\cong \F[x_1,\dots , x_n]^{\Sym_n}/(P_{d+1}, \dots , P_{d+n})$$
%where the $P_{d+1}, \dots , P_{d+n}$ form a regular sequence
%of maximal length in the ring of symmetric polynomials.
%This means that the ring of symmetric functions is in fact
%free and finitely generated over the polynomial subalgebra
%which they generate.
Thus the spectral sequence simplifies to
$$E_2^{*,*} =
Tor_{H^*(B,\F)}(\F, H^*(B,\F)
\otimes_{\F[u_1, \dots , u_n]}\F)
\cong
Tor_{\F[u_1, \dots , u_n]}(\F,\F)$$
This can be computed using
the standard Koszul complex, yielding
$$E_2 = \Lambda_\F(e_1, \dots, e_n)$$
where the $e_i$ are exterior classes in degree
$2r_i-1$.
There are no further differentials, as the algebra generators
for $E_2^{*,*}$ represent non--trivial elements in
the cohomology of $F$ which by construction
must transgress to the regular sequence
$\{u_1,\dots , u_n\}$ in $H^*(B, \F)$
in the Serre spectral sequence for the fibration
\[ F\to E \to B . \]
Therefore
the Eilenberg--Moore spectral sequence collapses at $E_2=E_{\infty}$.
Now this algebra is a free graded commutative algebra, hence
there are no extension problems and it follows that
$$H^*(F,\F)\cong \Lambda_\F (e_1, \dots , e_n)$$
as stated in the theorem.
\end{proof}

We now apply this result to the spaces $X(n,d)$.

\begin{thm}\label{thm0}
The cohomology of $X(n,d)$ is an exterior algebra on $n$ generators
$$H^*(X(n,d), \F)\cong
\Lambda_\F (E_{d+1}, \dots, E_{d+n}) \, , $$
where $E_j$ is a
cohomology class in dimension $2j-1$.
\end{thm}
\begin{proof}
As observed previously we have a fibration
\[
 X(n,d)\to (\mathbb CP^d)^n_{h\Sym_n}
\to BU(n) \, .
\]
The Eilenberg--Moore spectral sequence can therefore
be applied to compute the cohomology of $X(n,d)$.
The map
$F(n,d): (\mathbb CP^d)^n_{h\Sym_n}\to BU(n)$
induces a surjection of algebras
$$H^*(BU(n),\F)
\to H^*((\mathbb CP^d)^n_{h\Sym_n},\F)\to 0$$
which can be identified with the natural map
$$\F [x_1, \dots , x_n]^{\Sym_n}
\to (\F [x_1, \dots , x_n]/(x_1^{d+1}, \dots
x_n^{d+1}))^{\Sym_n}.$$
The kernel of this map is precisely the ideal
$$I_{n, d}= (x_1^{d+1}, \dots , x_n^{d+1})\cap  
\F [x_1, \dots , x_n]^{\Sym_n}.$$
By Theorem~\ref{thm2}(a), $I_{n, d}$ is generated by
a regular sequence of elements $P_{d+1},\dots ,P_{d+n}$. Here
each $P_j$ is a homogeneous polynomial in $x_1, \dots, x_n$ 
of degree $j$; its degree as
a cohomology class is $2j$.  These classes form a regular
sequence of maximal length in the polynomial
algebra $H^*(BU(n),\F)$. Thus the hypotheses of
Theorem~\ref{EMSS} hold, and the proof is complete.
\end{proof}

\begin{cor}
If $d<\infty$, then
$X(n,d)$ is a compact, connected, orientable
manifold.
\end{cor}
\begin{proof}
According to our calculation, for $m=n(n+2d)$ we
have $H^m(X(n,d),\mathbb Q)\cong\mathbb Q$. This
is precisely the dimension of the compact
manifold $X(n,d)=U(n)\times_{NT} (\mathbb S^{2d+1})^n$,
whence the result follows.
\end{proof}

\begin{remark}
Note that as $d$ gets large, the connectivity of the
space $X(n,d)$ increases; this is consistent with the stable
calculation, namely the acyclicity of $U(n)/NT$.
Also note that the manifold $U(n)/NT$ is not orientable,
as it is $\mathbb Q$--acyclic.
\end{remark}

For the case of $W(n,d)$ we offer the following general
result:

\begin{thm} \label{thm.serre}
For any field $\F$ of coefficients, the Serre spectral sequence for
the fibration $(\mathbb S^{2d+1})^n\to W(n,d)\to U(n)/T$
collapses at $E_2$ if and only if $d\ge n-1$, from which
we obtain an additive calculation
$$H^*(W(n,d),\F)\cong H^*(U(n)/T)\otimes H^*((\mathbb S^{2d+1})^n,\F).$$
In particular if $n!$ is invertible in $\F$, then
$$H^*(X(n,d),\F)\cong
[H^*(U(n)/T)\otimes H^*((\mathbb S^{2d+1})^n,\F)]^{\Sym_n}
\cong
\Lambda_\F (E_{d+1}, \dots, E_{d+n}) \, .$$
\end{thm}

\begin{proof}
Consider the Serre spectral sequence with $\F$
coefficients for the fibration
$(\mathbb S^{2d+1})^n\to W(n,d)\to U(n)/T$. The base is
simply connected and the cohomology of the fiber is
generated by the natural generators for the
$2d+1$--dimensional cohomology of each sphere.
The first differential in the spectral sequence can
be computed as follows: if $e_i\in
H^{2d+1}((\mathbb S^{2d+1})^n)$ is a natural generator then
$$d_{2d+2}(e_i) = [x_i^{d+1}] \in H^*(U(n)/T)\cong
H^*(BT)/(s_1, s_2,\dots ,s_{n})$$
where the $s_1, s_2,\dots ,s_{n}$ are the symmetric
polynomials. This follows from the diagram of fibrations
in the previous section
and the well-known calculation of the cohomology
of $(\mathbb CP^d)^n$ and $U(n)/T$ as quotients of
$H^*(BT)$. We now need the following algebraic lemma.

\begin{lem}\label{powers}
Let $\mathbb F$ be a commutative ring and $I$
be the ideal of $\mathbb F[x_1, ..., x_n]$ generated by the elementary
symmetric polynomials $s_1, ..., s_n$ in $x_1, ..., x_n$.
Then (a) $x_1^n \in I$ but (b) $x_1^{n-1} \not \in I$.
\end{lem}

Suppose Lemma~\ref{powers} is established (we only need it 
in the case where $\F$ is a field).
Then we conclude that $d_{2d+2}(e_i)=[x_i^{d+1}]=0$ in $H^*(U(n)/T)$
for all $i=1,\dots, n$ if and only if $d\ge n-1$. This 
implies that all the differentials in the spectral sequence are zero
and so it collapses at $E_2$. The assertions of Theorem~\ref{thm.serre}
follow from this and Theorem~\ref{thm0}.

\smallskip
It thus remains to prove Lemma~\ref{powers}.

(a) Recall that $x_1, \ldots , x_n$ are, by definition,
the roots of the polynomial
\[ x^n - x^{n-1} s_1 + x^{n-2} s_2 -  ... + (-1)^n s_n = 0 \, . \]
Thus $x_1^n = x_1^{n-1} s_1 - x^{n-2} s_2 + ... - (-1)^n s_n$,
and since every term in the right hand side lies in $I$,
part (a) follows.

\smallskip
(b) Assume, to the contrary, that
\begin{equation} \label{e.hypothetical-formula}
x_1^{n-1} = f_1 s_1 + ... + f_n s_n
\end{equation}
for some polynomials $f_1, ...,  f_n \in \mathbb F[x_1, \dots, x_n]$.
If such an identity is possible over $\mathbb F$, and 
$\alpha \colon \mathbb F \to L$ is a ring homomorphism
then, applying $\alpha$ to each of the coefficients of $f_1, \dots, f_n$,
we obtain an identity of the same form over $L$. Thus, for the purpose
of showing that~\eqref{e.hypothetical-formula} is not possible,  
we may, without loss of generality, replace $\mathbb F$ by $L$.
In particular, we may take $L$ to be the algebraic closure of
the field $\mathbb F/M$, where $M$ is a maximal 
ideal of $\mathbb F$. After replacing $\mathbb F$ by this $L$, 
we may assume that $\mathbb F$ is an algebraically closed field.

Equating the homogeneous terms of degree $n-1$ on both sides, we
see that after replacing $f_1, f_2, \dots, f_{n-1}$ by their
homogeneous parts of degrees $n-2, n-3, \ldots, 0$, respectively,
we may assume that $f_n = 0$.

Since $\mathbb F$ is an algebraically closed field, 
$x^n - 1$ factors into a product of linear terms 
\begin{equation} \label{e.zeta}
x^n - 1 = (x - \zeta_1)(x - \zeta_2) \cdot \ldots \cdot (x - \zeta_n)\, .
\end{equation}
for some $\zeta_1, \dots, \zeta_n \in \mathbb F$.
(As an aside, we remark that $\zeta_1, \dots, \zeta_n \in \mathbb F$ 
are distinct
if $p = \Char(\mathbb F)$ does not divide $n$ but not in general;
at the other extreme, if $n$ is a power
of $p$ then $\zeta_1 = \dots = \zeta_n = 1$.) By~\eqref{e.zeta}
\[ \text{$s_i(\zeta_1, \dots, \zeta_n) = (-1)^i$ (coefficient of $x^{n-i}$
in $x^n - 1$) = 0} \]
for every $i = 1, \dots, n-1$.
Hence, substituting $\zeta_i$ for $x_i$ in~\eqref{e.hypothetical-formula}, 
and remembering that $f_n = 0$,
we obtain $\zeta_1^{n-1} = 0$, i.e., $\zeta_1 = 0$. Since $\zeta_1$ is a root
of $x^n - 1 = 0$, we have arrived at a contradiction.
This shows that~\eqref{e.hypothetical-formula} is impossible.
The proof of Lemma~\ref{powers} and thus of Theorem~\ref{thm.serre}
is now complete.
\end{proof}

Calculations with field coefficients can be pieced
together to provide information on the integral
cohomology of $X(n,d)$.

\begin{prop}
The cohomology ring $H^*(X(n,d),\mathbb Z)$ has no $p$--torsion if
$p>n$.
\end{prop}
\begin{proof}
By our previous results if $p>n$ then
$$\rm{dim}_{\mathbb F_p}~H^*(X(n,d),\mathbb F_p)
= \rm{dim}_{\mathbb Q}~H^*(X(n,d),\mathbb Q)=2^n.$$ Hence by the
universal coefficient theorem, there can be no $p$--torsion in the
integral cohomology of $X(n,d)$.
\end{proof}

The situation is more complicated if $n \ge p = \Char(\F)$.
In particular, we will show that in this case
the kernel $I(n, d)$ of the map $H^*(BU(p),\F_p)
\to H^*((\mathbb CP^d)^p,\F_p)$ cannot be generated by
a regular sequence for any $d \ge 2$ (and, in most cases
for $d = 1$ as well); see Theorem~\ref{thm2}(b).
%%%%%%%%%%%%%%%%%%%%%%%%%%%%%%%%
% \begin{prop}
% The kernel $I(n,d)$ of the map $H^*(BU(p),\F_p)
% \to H^*((\mathbb CP^d)^p,\F_p)$ cannot be generated by
% a regular sequence. Moreover a minimal set of generators
% has $p+1$ elements.
% \qed
% \end{prop}
%%%%%%%%%%%%%%%%%%%%%%%%%%%%%%%%
We now provide an explicit calculation in the case where $n = d = p = 2$.

\begin{example}
Consider the map
$\tilde{F}(2,2): \mathbb S^2\times\mathbb S^2 \to BU(2)$.
Its fiber is 
\[ W(2,2) = U(2)\times_T (\mathbb S^3\times\mathbb
S^3) \]
which itself fibers over $U(2)/T=\mathbb S^2$ with fiber $\mathbb
S^3\times\mathbb S^3$. Hence for dimensional reasons
$H^*(W(2,2),\mathbb Z)\cong H^*(\mathbb S^3\times\mathbb S^3\times
\mathbb S^2,\mathbb Z)$. The $S_2$--action on this space exchanges
the two $3$-spheres and applies the antipodal map on $\mathbb S^2$.
Thus the orbit space $X(2,2)$ will be rationally cohomologous to
$\mathbb S^3\times\mathbb S^5$, as predicted by Theorem~\ref{thm0}.
However, it can be shown that $H^*(X(2,2),\F_2)$ has Poincar\'e
series
$$p(t) = 1 + t + t^2+t^3+t^5+t^6+t^7+t^8 \, . $$
On the other hand, the corresponding Poincar\'e series for rational
cohomology is
$$q(t) = 1+ t^3 + t^5 + t^8$$
which accounts for the torsion free classes in the integral
cohomology. This example illustrates the presence of
$2$--torsion in the cohomology of $X(2,2)$. Of course in
this case we have $\pi_1(X(2,2))=\mathbb Z/2$, which accounts
for the classes in degrees one and two in mod 2 cohomology,
and by Poincar\'e duality for the classes in degrees six
and seven.

On the other hand, recall that if $H^*(BU(2),\F_2)\cong \F_2[c_2, c_4]$
and $H^*(\mathbb S^2\times\mathbb S^2,\F_2)\cong
\Lambda (u_2, v_2)$ then
$\tilde{F}(2,2)^*(c_2)=u_2+v_2$ and $\tilde{F}(2,2)^*(c_4)=u_2v_2$. 
Thus we see that $\tilde{F}(2,2)^*$
is not surjective and that its kernel is
generated by the classes
$c_2^2, c_2^3 + c_2c_4, c_4^2$. These classes correspond
to the symmetric polynomials
$P_2 = x_1^2+x_2^2$, $P_3 = x_1^3+x_2^3$ and $P_{2, 2} = x_1^2x_2^2$.
Note that if $2$ is invertible in the coefficients then
$$P_{2,2} = \frac{P_2^2-(x_1+x_2)P_3 + (x_1x_2)P_2}{2} \, , $$
and the third generator is redundant.
\end{example}

More generally, using the algebraic calculations in
Theorem~\ref{thm1}, Theorem~\ref{thm2} and Corollary~\ref{cor.thm2c} 
we obtain the following.

\begin{thm}
Assume that $p\le n \le 2p -1$ and $d\ge 2$. Then
the kernel of the map induced by $\tilde{F}(n,d)$ in cohomology
$$\tilde{F}(n,d)^*: H^*(BU(n),\F_p)\to H^*((\mathbb CP^d)^n,\F_p)$$
is generated by the following $n + 1$ elements:
\begin{itemize}
\item $P_{d+i}$, where $1\le i\le n$ and $|P_j| = 2j$
\item $P_{\underbrace{\text{\tiny $d+1, \ldots, d+1$}}_{\text{\tiny $p$ times}}}$ and $|P_{\underbrace{\text{\tiny $d+1, \ldots,
d+1$}}_{\text{\tiny $p$ times}}}| = 2p(d+1)$
\end{itemize}
Moreover these elements cannot form a regular sequence, and this kernel
cannot be generated by less than $n+1$ elements.
\qed
\end{thm}

\section{The orthogonal groups and more calculations at $p=2$}

The situation for $p=2$ is somewhat different, as there are
specific geometric models which
are special to this characteristic. Here we consider
the standard diagonal inclusion
$V=(\mathbb Z/2)^n \hookrightarrow O(n)$
into the group of orthogonal $n\times n$ matrices.
The group $V$ is self-centralizing in $O(n)$; its normalizer
$NV$ is the wreath product $NV = \mathbb Z/2 \wr \Sym_n$.
The Weyl group $W = NV/V$
of $V$ in $O(n)$ is thus isomorphic to $\Sym_n$; it acts on
$V=(\mathbb Z/2)^n$ by permuting the $n$ factors of $\mathbb Z/2$.
The classifying space for $V$ is
$BV = (\mathbb RP^\infty)^n$, its mod $2$
cohomology is a polynomial algebra on $n$ one dimensional
generators $\mathbb F_2[x_1,\dots , x_n]$. The inclusion
induces a map from the cohomology of $BO(n)$ to this
algebra,
which gives rise to an isomorphism
onto the symmetric invariants.
As before, the truncated projective space $\mathbb RP^d$ is a
natural subspace of $\mathbb RP^\infty$, and Theorem~\ref{thm1}
provides a description of the kernel of the homomorphism induced by
the map $H(n,d): (\mathbb RP^d)^n\to BO(n)$ for $n=1,2,3$.

The classifying space for $NV=\mathbb Z/2 \, \wr \,  \Sym_n$ is
$BNV = (\mathbb RP^\infty)^n_{h\Sym_n}$. However, 
as our calculations are at $p=2$ and $|\Sym_n|$ is even, 
the homotopy orbit space has a lot more cohomology than just
the truncated symmetric invariants (for example, it contains
a copy of $H^*(\Sym_n,\F_2)$).
The wreath product $NV$ acts on $(\mathbb S^d)^n$ extending
the coordinatewise antipodal action of $V$.
Thus we have a fiber bundle
$(\mathbb S^d)^n \to
(\mathbb RP^d)^n_{h\Sym_n}
\to BNV$
where we identify  $(\mathbb S^d)^n_{hNV}
\simeq (\mathbb RP^d)^n_{h\Sym_n}$.

\begin{example}
For $n=2$ we can identify $NV$ with the dihedral 
group $D_8$ and 
its cohomology has generators
$e$, $u$  $v$ in degrees $1,1,2$ respectively with
the single relation $e\cdot u =0$ (see \cite{AM}). 
The elements $u,v$
can be identified with the standard symmetric generators 
$x_1+x_2$ and $x_1x_2$
in $H^*(V,\F_2)^{\Sym_2}$ via the restriction map. 
In fact we
have isomorphisms (see \cite{AM}, page 118)
$H^*(BD_8,\mathbb F_2)\cong H^*(\Sym_2,H^*(V,\mathbb F_2))$
and
$H^*((\mathbb S^d)^2_{hD_8}, \F_2)\cong H^*(\Sym_2,
H^*((\mathbb RP^d)^2, \F_2)$.
Using these descriptions and
Theorem~\ref{thm1} it can be shown that the 
the kernel of the homomorphism
$H^*(BD_8, \F_2)\to H^*((\mathbb S^d)^2_{hD_8}, \F_2)$
is the
ideal
generated by the three elements $P_{d+1}=x_1^{d+1}+x_2^{d+2}$,
$P_{d+2}=x_1^{d+2}+x_2^{d+2}$ and $P_{d+1,d+1}=x_1^{d+1}x_2^{d+1}$.
This ideal is called the Fadell--Husseini index (see \cite{FH})
of the $D_8$--space
$\mathbb S^d\times\mathbb S^d$;
it has some interesting applications in topology
and it has been fully calculated
in \cite{BZ}. 
\end{example}

Geometrically, the fibration which our mod $2$ calculations can be 
applied to is described by the diagram:
\[
\xymatrix{
& (\mathbb S^d)^n\ar@{=}[r] \ar[d] & (\mathbb S^d)^n \ar[d] \\
Y(n,d)\ar@{=}[r] & O(n)\times_V(\mathbb S^d)^n \ar[r] \ar[d] &
(\mathbb RP^d)^n \ar[r]^{H(n,d)} \ar[d] & BO(n)\ar@{=}[d] \\
& O(n)/V\ar[r]  & BV \ar[r]
& BO(n)\\
}
\]

Here we recall some classical results. First, from the homotopy
long exact sequence of the fibration we see that $O(n)/V$ is
path--connected because $\pi_1(BV)\to \pi_1(BO(n))\cong\mathbb Z/2$ is 
surjective
(the dual map in mod $2$ cohomology is injective). 
Its fundamental group
acts homologically trivially on $H^*((\mathbb S^d)^n ,\F_2)$,
as it acts through its image in $V$. Therefore the Serre spectral
sequence for the fibration
$(\mathbb S^d)^n \to Y(n,d)
\to O(n)/V$ 
has the form
$$E_2^{*,*} = H^*(O(n)/V)\otimes H^*((\mathbb S^d)^n,\F_2)
\implies H^*(Y(n,d),\F_2) \, . $$
Using Lemma~\ref{powers}, we see that this spectral sequence collapses 
at $E_2$ if and only if $d\ge n-1$. 

\begin{thm}
If $d\ge n-1$ then
we have an additive
isomorphism
$$H^*(Y(n,d),\F_2)\cong H^*(O(n)/V)\otimes H^*((\mathbb S^d)^n,\F_2).$$
\qed
\end{thm}

\smallskip
\section{Truncated symmetric polynomials}
\label{sect5}
We will now state and prove the algebraic results used 
in the previous sections. For the rest of the paper
we will use the following notations:
\[ \begin{array}{lcl}
n, d          & & \text{positive integers} \\
\Sym_n          & & \text{the symmetric group on $n$ letters} \\
\mathbb F          & & \text{base field} \\
x_1, \dots, x_n  & & \text{independent variables over $\mathbb F$} \\
\mathbb F[x_1, \dots, x_n]  & & \text{polynomial ring in $n$ variables} \\
R_n := \mathbb F[x_1, \dots, x_n]^{\Sym_n}  & & 
\text{ring of symmetric polynomials in $n$ variables} \\
I_{n, d} := (x_1^{d + 1}, \dots, x_n^{d + 1}) \cap R_n & & 
\text{ideal of $R_n$.} 
\end{array} \]

If $a_1, \dots, a_n$ are non-negative integers, we will write
$P_{a_1, \dots, a_n}$ for the orbit sum of $x_1^{a_1} \dots x_n^{a_n}$.
In other words, $P_{a_1, \dots, a_n}$ is the sum of
monomials $x^{a_1'} \dots x_n^{a_n'}$, as
$a_1', \dots, a_n'$ range over all possible permutations of
$a_1, \dots, a_n$. This sum
has $\dfrac{n!}{\lambda_1!  \cdots \lambda_m!}$ terms,
where $\lambda_1, \dots, \lambda_m$ is the partition of $n$ associated to
$a_1, \dots, a_n$. (Recall that this means that that there are $m$
distinct integers among $a_1, \dots, a_n$, occurring with multiplicities
$\lambda_1, \dots, \lambda_m$, respectively.)

Permuting $a_1, \dots, a_n$ does not change 
$P_{a_1, \dots, a_n}$, so we will always assume
that $a_1 \ge \dots \ge a_n$. With this convention, 
the orbit sums $P_{a_1, \dots, a_n}$
clearly form a basis of $R_n := \mathbb F[x_1, \dots, x_n]^{\Sym_n}$
as an $\mathbb F$-module. The multiplication rule in this basis 
is given by
\begin{equation} \label{e-1}
P_{a_1, \dots, a_n} \cdot P_{b_1, \dots, b_n} =
\sum k_{c_1, \dots, c_n} P_{c_1, \dots, c_n} \, ,
\end{equation}
where $c_1 \ge \ldots \ge c_n$ and there 
are exactly $k_{c_1, \dots, c_n}$ different ways to write
\[ (c_1, \dots, c_n) = (a_1', \dots, a_n') + (b_1', \dots, b_n') \]
for some permutation
$a_1', \dots, a_n'$ of $a_1, \dots, a_n$ and
some permutation $b_1', \dots, b_n'$ of $b_1, \dots, b_n$.

To make our formulas less cumbersome, we will often
abbreviate $P_{a_1, \dots, a_r, 0, \dots, 0}$
as $P_{a_1, \dots, a_r}$. As long as the number of variables $n$ is
fixed, this will not lead to any confusion. For example, in this 
notation,
\[ P_i = x_1^i + \dots + x_n^i \]
is the usual power sum of degree $i$ and
\begin{equation} \label{e.esp}
\begin{array}{l} P_1 = x_1 + \dots + x_n, \\
P_{1, 1} = x_1 x_2 + \dots + x_{n-1}x_n, \\
\dots \\
P_{\underbrace{\text{\tiny $1, \ldots, 1$}}_{\text{\tiny $n$ times}}} =
x_1 x_2 \dots x_n
\end{array} 
\end{equation}
are the elementary symmetric polynomials.

% Our main algebraic results are the following two theorem.

The main result of this section is the following theorem.

\begin{thm} \label{thm1} Let $\F$ be a field of characteristic $p \ge 0$.

\smallskip
(a) If $p = 0$ or $n < p$ then the ideal
$I_{n, d} := (x_1^{d+1}, \dots, x_n^{d+1}) \cap 
\mathbb F[x_1, \dots, x_n]^{\Sym_n}$ 
of $R_n := \mathbb F[x_1, \dots, x_n]^{\Sym_n}$ 
is generated by $P_{d +1}, \dots, P_{d + n}$.

\smallskip
(b) If $n \le 2p -1$ then $I_{n, d}$
is generated by $P_{d +1}, \dots, P_{d + n}$ and
$P_{\underbrace{\text{\tiny $d+1, \ldots, d+1$}}_{\text{$p$ times}}}$.
\end{thm}

The rest of this section will be devoted to proving Theorem~\ref{thm1}.
First we note that every element of $I_{n, d}$ is an
$\mathbb F$-linear combination of orbit sums $P_{a_1, \dots, a_n}$,
where $a_1 \ge d+1$. Thus in order to prove 
Theorem~\ref{thm1} it suffices to show that every
$P_{a_1, \dots, a_n}$ lies in $I$. Our first step in
this direction is the following lemma.

We define the {\em weight} of the orbit sum $P_{a_1, \dots, a_n}$ 
as the largest integer $r \le n$ such that $a_r \ge 1$.
As mentioned above, we will abbreviate such an orbit sum
as $P_{a_1, \dots, a_r}$.

We define the {\em leading multiplicity} of $P_{a_1, \dots, a_n}$ 
as the largest integer $s \le n$ such that $a_1 = \dots = a_s$.
Here, as always, we are assuming that $a_1 \ge a_2 \ge \ldots \ge a_n \ge 0$.

\begin{lem} \label{lem3} Let $\F$ be a field and $J_{n, d}$ be the ideal of 
$R_n = \mathbb F[x_1, \dots, x_n]^{\Sym_n}$ generated 
by $P_{d + 1}, \dots, P_{d + n}$. Then 
$J_{n, d}$ contains every orbit sum  $P_{a_1, \dots, a_n}$ 
with $a_1 \ge d + 1$, whose leading multiplicity is invertible in $\F$.
\end{lem}

The leading multiplicity of $P_{a_1, \dots, a_n}$
is, by definition, an integer between 
$1$ and $n$.  Theorem~\ref{thm1}(a) is thus an immediate consequence 
of this lemma.

\begin{proof} We will argue by induction on the weight
$r$ of $P_{a_1, \dots, a_n}$. For the base case, let $r = 1$.
That is, we claim that
$P_i \in J_{n, d}$ for every $i \ge d + 1$. For $i = d+1, \dots, d+n$ 
this is given.  Applying Newton's identities,
\[ P_{m + n + 1} = P_1 \cdot P_{m + n} - P_{1, 1}  \cdot P_{m + n-1} +
\dots + (-1)^{n+1} 
P_{\underbrace{\text{\tiny $1, \ldots, 1$}}_{\text{$n$ times}}} 
\cdot P_{m + 1} \]
recursively, with $m = d, d+ 1, d+2$, etc., we see that
$P_{m + n +1} \in J_{n, d}$ for every $m \ge d$. This settles the base
case.

For the induction step assume that $r \ge 2$.
By~\eqref{e-1},
\begin{equation} \label{e1}
P_{a_1} \cdot P_{a_2, \dots, a_r} =
s P_{a_1, a_2, \dots, a_r} + 
P_{a_1 + a_2, a_3, \dots, a_r} + 
P_{a_1 + a_3, a_2, a_4,  \dots, a_r} + \dots +
P_{a_1 + a_r, a_2, a_3,  \dots, a_{r-1}} \, . 
\end{equation}
Each of the terms
\[ P_{a_1 + a_2, a_3, \dots, a_r}, \;
P_{a_1 + a_3, a_2, a_4,  \dots, a_r}, \; \ldots, \;
P_{a_1 + a_r, a_2, a_3,  \dots, a_{r-1}}  \]
is an orbit sum of leading multiplicity $1$
and weight $r - 1$. By the induction 
assumption each of them lies in $J_{n, d}$.
Since we also know that $P_{a_1} \in J_{n, d}$,
equation~\eqref{e1} tells us that $P_{a_1, \dots, a_r} 
\in J_{n, d}$ whenever $s$ is invertible in $\F$.
\end{proof}

We are now turn to the proof of Theorem~\ref{thm1}(b).
In view of part (a), we may assume that $p < n \le 2p - 1$.
Let $I$ be the ideal of $R_n = \mathbb F[x_1, \dots, x_n]^{\Sym_n}$ 
generated by the polynomials listed in the statement of 
Theorem~\ref{thm1}(b).  Recall that it suffices to show that
\begin{equation} \label{e.goal}
\text{$P_{a_1, \dots, a_n} \in I$ whenever $a_1 \ge d + 1$.}
\end{equation} 
Denote the leading multiplicity of
$P_{a_1, \dots, a_n}$ by $s$. We will now consider three cases.

\smallskip
{\bf Case 1.} $s \ne p$. 
Since we are assuming that $p < n \le 2p - 1$, this is equivalent to
$s$ being invertible in $\F$.   Clearly, $J_{n, d} \subset I$;
Lemma~\ref{lem3} thus tells us that~\eqref{e.goal} holds.

\smallskip
{\bf Case 2.} $s = p$ and $P_{a_1, \dots, a_n}$ has weight $p$. In other 
% words, $(a_1, \dots, a_n) = (a, \dots, a, 0, \dots, 0)$, 
% where $a \ge d+1$ is repeated $p$ times. We 
words, we want to show that 
\begin{equation} \label{e2.0}
P_{\underbrace{\text{\tiny $a, \dots, a$}}_{\text{\tiny $p$ times}}} 
\in I \, . 
\end{equation}
Let $e = a - (d + 1)$. By~\eqref{e-1} we see that
\begin{equation} \label{e2}
P_{\underbrace{\text{\tiny $d+1, \dots, d + 1$}}_{\text{\tiny $p$ times}}} 
\cdot
P_{\underbrace{\text{\tiny $e, \dots, e$}}_{\text{\tiny $p$ times}}} =
P_{\underbrace{\text{\tiny $a, \dots, a$}}_{\text{\tiny $p$ times}}} +
\Gamma \, , \end{equation}
where $\Gamma$ is a positive integer linear combination of orbit sums
of leading multiplicity $\le p-1$. By Lemma~\ref{lem3} 
$\Gamma \in I$. Since by definition, 
$P_{\underbrace{\text{\tiny $d+1, \dots, d + 1$}}_{\text{\tiny $p$ times}}}$
lies in $I$, the left hand side also lies in $I$. Thus~\eqref{e2.0}
holds as well.

Note that the above argument depends, in a crucial way, on
our assumption that $n \le 2p - 1$. For $n \ge 2p$ the sum $\Gamma$
in~\eqref{e2} would contain a term of the form 
$P_{d+1, \dots, d+1, e, \dots, e}$ (or $P_{e, \dots, e, d+1, \dots, d+1}$, 
if $e > d + 1$), with each $e$ and $d+1$ repeating exactly $p$ times.
This orbit sum has leading multiplicity $p$, 
and in the case we cannot conclude that $\Gamma \in I$.

\smallskip
{\bf Case 3.} $s = p$, general case. Denote $a_1 = \dots = a_p$ by $a$.
Using formula~\eqref{e-1} once again, we see that
\[ P_{a_1, \dots, a_n} =
P_{\underbrace{\text{\tiny $a, \dots, a$}}_{\text{\tiny $p$ times}}}
\cdot P_{a_{p+1}, \dots, a_n} + \Delta \, , \]
where $\Delta$ is an integer linear combination of orbit
sums $P_{c_1, \dots, c_n}$ of leading multiplicity $\le p -1$.
Note that
$P_{\underbrace{\text{\tiny $a, \dots, a$}}_{\text{\tiny $p$ times}}}
\in I$ by Case 2 and $\Delta \in I$ by Lemma~\ref{lem3}.
We thus conclude that $P_{a_1, \dots, a_p} \in I$ as well. This completes
the proof of Theorem~\ref{thm1}.
\qed

\section{Regular sequences} 
\label{sect.reg-seq}

We now turn to the question of whether or not the ideal 
$I_{n, d} = (x_1^{d+1}, \dots, x_n^{d+1}) \cap R_n$ of
$R_n = \mathbb{F}[x_1, \dots, x_n]^{\Sym_n}$
can be generated by a regular sequence. Our goal is to prove 
the following theorem.

\begin{thm} \label{thm2} Let $\mathbb F$ be a field of characteristic
$p \ge 0$.

\smallskip
(a) If $n!$ is not divisible by $p$ then $I_{n, d}$ is generated
by the regular sequence $P_{d + 1}, \dots, P_{d + n}$ in $R_n$.

\smallskip
(b) Assume that $n \ge p > 0$ and either (i) $n \not \equiv -1$ (mod $p$)
and $d \ge 1$ or (ii) $n \equiv -1$ (mod $p$) and $d \ge 2$.
Then $I_{n, d}$ is not generated by any regular sequence in $R_n$.
\end{thm}

The assumptions on $d$ in part (b) cannot be dropped;
see Remark~\ref{rem.thm2}.  Our proof of Theorem~\ref{thm2} will rely on
the following elementary lemma.

\begin{lem} \label{lem.thm2}
(a) The elements $P_{a_1, \dots, a_n}$, with
$d \ge a_1 \ge \dots \ge a_n \ge 0$ form a basis for 
$R_n/I_{n, d}$ as an $\mathbb F$-vector space.   

\smallskip
(b) The Krull dimension of $R_n/I_{n, d}$ is $0$.

\smallskip
(c) Suppose $I_{n, d}$ is generated by $r_1, \dots, r_m \in R_n$,
as an ideal of $R_n$. Then $m \ge n$. Moreover, $r_1, \dots, r_m$ form 
a regular sequence in $R_n$ if and only if $m = n$. 
\end{lem}

\begin{proof} (a) The power sums $P_{a_1, \dots, a_n}$ with
$a_1 \ge \dots \ge a_n \ge 0$ form an $\mathbb F$-basis of $R_n$. 
The power sums $P_{a_1, \dots, a_n}$ with
$a_1 \ge \dots \ge a_n \ge 0$ and $a_1 \ge d + 1$  
form an $\mathbb F$-basis of $I_{n, d}$, and part (a) follows. 

(b) By part (a) $R_n/I_{n, d}$ is a finite-dimensional 
$\mathbb F$-vector space.

(c) Recall that $R_n$ is a polynomial
ring over $\mathbb F$ generated by the elementary symmetric polynomials
in $x_1, \dots, x_n$. In particular, $R_n$ is a Cohen-Macauley ring. 
Part (c) now follows from part (b).  
\end{proof}

\begin{proof}[Proof of Theorem~\ref{thm2}]
(a) If $p = \Char(\mathbb F)$ does not divide $n!$ then Theorem~\ref{thm1}(a)
tells us that $I_{n, d}$ is generated, as an ideal of $R_n$, 
by the $n$ elements $P_{d + 1}, \dots, P_{d + n}$. 
By Lemma~\ref{lem.thm2}(c) these elements form 
a regular sequence in $R_n$.  

\smallskip
(b) If $I_{n, d}$ is generated by a regular 
sequence then $\Soc(R_n/I_{n, d})$ is a 1-dimensional
$\mathbb F$-vector space; see, e.g.~\cite[p. 144]{ns} 
or~\cite[Section 21.2]{eisenbud}. It is an immediate 
consequence of the multiplication formula~\eqref{e-1} that
\[ P_{\underbrace{\text{\tiny $d, \dots, d$}}_{\text{\tiny $n$ times}}} \in
\Soc(R_n/I_{n, d}) \] for any $\mathbb F$, $d$ and $n$. Thus 
in order to show that $I_{n, d}$ is not generated by a regular 
sequence it suffices to exhibit an orbit sum
$P_{a_1, \dots, a_n} \in \Soc(R_n/I_{n, d})$, with
$(a_1, \dots, a_n) \ne (d, \dots, d)$. (Indeed,
$P_{a_1, \dots, a_n}$ and 
$P_{\underbrace{\text{\tiny $d, \dots, d$}}_{\text{\tiny $n$ times}}}$
are $\mathbb F$-linearly independent in $R_n/I_{n, d}$ by
Lemma~\ref{lem.thm2}(a).)

\smallskip
(i) Suppose $d \ge 1$ and $n  = pq + r$, where $q \ge 1$ and
$r \in \{ 0, 1, \dots, p - 2 \}$.  
We claim that in this case $P_{a_1, \dots, a_n}$
lies in $\Soc(R_n/I_{n, d})$, if 
\[ \text{$a_1 = \dots = a_{pq - 1} = d$
and $a_{pq} = a_{pq+1} = \dots = a_n = d-1$.} \]

In other words, we claim that
%%%%%%%%%%%%%%%%%%%%%%%%
% $P_{\underbrace{\text{\tiny $d, \dots, d,$}}_{\text{\tiny $pq-2$ times}}
% {\text{\tiny $d-1, \dots, d-1$}}_{\text{\tiny $r+1$ times}}}$
% lies in $\Soc(R_n/I_{n, d})$. In other words,
% \[ P_{\underbrace{\text{\tiny $d, \dots, d$}}_{\text{\tiny $pq-2$ times}}
% \cdot P_{a_1, \dots, a_n} \in I_{n, d} \]
%%%%%%%%%%%%%%%%%%%%%%%%
\begin{equation} \label{e.socle1}
P_{a_1, \dots, a_n} \cdot P_{b_1, \dots, b_n} \in I_{n, d} 
\end{equation}
for every $b_1 \ge b_2 \ge \dots \ge b_n \ge 0$, with $b_1 \ge 1$.
To prove this claim, we will examine the product 
$P_{a_1, \dots, a_n} \cdot P_{b_1, \dots, b_n}$
using the multiplication formula~\eqref{e-1}.

First of all, note that we may assume without loss of generality that 
\begin{equation} \label{e.socle1.5}
(b_1, \dots, b_n) =
({\underbrace{1, \dots, 1}_{\text{\tiny $s$ times}}}, 0, \dots, 0)
\end{equation}
for some $1 \le s \le r + 1$. For all other choices of 
$(b_1, \dots, b_n) \ne (0, \dots, 0)$ every term $P_{c_1, \dots, c_n}$ 
appearing in the right hand side of the formula~\eqref{e-1} 
will have $c_1 \ge d+1$ and thus will lie in $I_{n, d}$ (for any
base field $\mathbb F$).

If $(b_1, \dots, b_n)$ is as in~\eqref{e.socle1.5}, 
the only orbit sum $P_{c_1, \dots, c_n}$ with $c_1 \le d$
appearing in the right hand side of~\eqref{e-1}, will 
have $c_1 = \dots = c_{pq + s - 1} = d$ and $c_{pq + s} = c_{pq + s + 1} = 
\dots = c_n = d-1$. This sum will appear with coefficient
$k_{c_1, \dots, c_n}$ = number of ways to write $(c_1, \dots, c_n)$
as $(a_1', \dots, a_n') + (b_1', \dots, b_n')$, where
$(a_1', \dots, a_n')$ is a permutation of $(a_1, \dots, a_n)$ and
$(b_1', \dots, b_n')$ is a permutation of $(b_1, \dots, b_n)$. 
We claim that $k_{c_1, \dots, c_n}$ is divisible by $p$
and hence, is $0$ in $\mathbb F$; this will immediately 
imply~\eqref{e.socle1}.  Indeed, in this case
$k_{c_1, \dots, c_n}$ is simply the number of ways to specify which
$s$ of the elements $b_1', \dots, b_{pq+s -1}'$ should be equal to $1$
(the remaining ones will be $0$). Thus
\[ k_{c_1, \dots, c_n} = {pq + s - 1 \choose s} \, . \] 
Since $q \ge 1$ and $1 \le s \le r+1 \le p-1$, this number 
is divisible by $p$, as claimed.

\smallskip
(ii) Now suppose $d \ge 2$ and $n  = pq + p - 1$, where $q \ge 1$.
We claim that in this case $P_{a_1, \dots, a_n}$
lies in $\Soc(R_n/I_{n, d})$, if $a_1 = \dots = a_{pq - 1} = d$
and $a_{pq} = a_{pq+1} = \dots = a_{pq + p - 2} = d-1$ and 
$a_{pq + p -1} = d - 2$. 

Once again, we need to show that~\eqref{e.socle1}
holds for every $b_1 \ge b_2 \ge \dots \ge b_n \ge 0$, with $b_1 \ge 1$.
The analysis of the product $P_{a_1, \dots, a_n} \cdot P_{b_1, \dots, b_n}$,
based on formula~\eqref{e-1}, is similar to part (i) but 
a bit more involved. 
First of all, we may assume without loss of generality that 
\begin{equation} \label{e.socle2}
(b_1, \dots, b_n) =
(2, {\underbrace{1, \dots, 1}_{\text{\tiny $t$ times}}}, 0, \dots, 0)
\end{equation}
for some $t \in [0, p-1]$ or
\begin{equation} \label{e.socle3}
(b_1, \dots, b_n) =
({\underbrace{1, \dots, 1}_{\text{\tiny $s$ times}}}, 0, \dots, 0)
\end{equation}
for some $s \in [1, p]$.
For other $(b_1, \dots, b_n) \ne (0, \dots, 0)$, every orbit sum
$P_{c_1, \dots, c_n}$ appearing in the right hand side of~\eqref{e-1}
will lie in $I_{n, d}$, so that~\eqref{e.socle1} will hold 
over any base field $\mathbb F$.

In case~\eqref{e.socle2} the only term $P_{c_1, \dots, c_n}$ with
$c_1 \le d$ appearing in the right hand side of~\eqref{e-1} will
have $c_1 = \dots = c_{pq + t} = d$ and $c_{pq + t +1} = \dots = c_{pq + p -1}
= d-1$. This term will appear with coefficient
$k_{c_1, \dots, c_n}$ = number of ways to write $(c_1, \dots, c_n)$
as $(a_1', \dots, a_n') + (b_1', \dots, b_n')$, where
$(a_1', \dots, a_n')$ is a permutation of $(a_1, \dots, a_n)$ and
$(b_1', \dots, b_n')$ is a permutation of $(b_1, \dots, b_n)$. 
One of the elements
$b_1', \dots, b_{pq + t}'$ should be equal to $2$, 
$t$ of these elements should be equal to $1$, and 
the remaining $pq - 1$ should be equal to $0$.
Thus
\[ k_{c_1, \dots, c_n} = (pq + t){pq + t - 1 \choose t} \, . \]
If $t = 0$ then the first factor is divisible by $p$. If $1 \le t \le p-1$
then the second factor is divisible by $p$.
Either way, $k_{c_1, \dots, c_n} = 0$ in $\mathbb F$, as desired.

In case~\eqref{e.socle3} exactly two orbit sums will appear in the 
right hand side of~\eqref{e-1}, namely
%  $P_{c_1, \dots, c_n}$ will
% appear in the right hand side of~\eqref{e-1}, with
\[ P_{\underbrace{\text{\tiny $d, \dots, d,$}}_{\text{\tiny $pq + s - 2$}}
\underbrace{\text{\tiny $d-1, \dots, d-1$}}_{\text{\tiny $p - s + 1$}}} \]
and
\[ P_{\underbrace{\text{\tiny $d, \dots, d,$}}_{\text{\tiny $pq - 1 + s$}}
\underbrace{\text{\tiny $d-1, \dots, d-1$}}_{\text{\tiny $p - s - 1$}}, 
\; d - 2}  \]
with coefficients
\[ k_{\underbrace{\text{\tiny $d, \dots, d,$}}_{\text{\tiny $pq+ s -2$}}
\underbrace{\text{\tiny $d-1, \dots, d-1$}}_{\text{\tiny $p - s + 1$}}} 
= {pq + s - 2 \choose s - 1} (p - s + 1) \]
and
\[ k_{\underbrace{\text{\tiny $d, \dots, d,$}}_{\text{\tiny $pq+ s -1$}}
\underbrace{\text{\tiny $d-1, \dots, d-1$}}_{\text{\tiny $p - s - 1$}}, \; 
d - 2} = {pq + s - 1 \choose s} \, , \]
respectively. (The second orbit sum does not occur if $s = p$.)
Both of these coefficients are divisible by $p$ and
hence, are $0$ in $\mathbb F$. This completes the proof 
of Theorem~\ref{thm2}.
\end{proof}

\begin{cor} \label{cor.thm2c}
Suppose (i) $p \le n \le 2p-2$ and $d \ge 1$ or (ii) $n = 2p - 1$ and $d \ge 2$.
Then the ideal $I_{n, d}$ can be generated by $n + 1$ elements of $R_n$ but
cannot be generated by $n$ elements.
\end{cor}

\begin{proof}
Theorem~\ref{thm1}(b) tells us that $I_{n, d}$ is generated
by $n + 1$ elements.
If $I_{n, d}$ could be generated by $n$ elements then
by Lemma~\ref{lem.thm2}(c) these $n$ elements 
would form a regular sequence in $R_n$, contradicting
Theorem~\ref{thm2}(b).
\end{proof}

\begin{remark} \label{rem.thm2} The conditions that $d \ge 1$ 
and $d \ge 2$ in parts (i) and (ii) of Theorem~\ref{thm2}(b) 
respectively, cannot be dropped. The same goes for conditions
(i) and (ii) in Corollary~\ref{cor.thm2c}.

Indeed, suppose $d = 0$. Recall that
$R_n = \mathbb F[x_1, \dots, x_n]^{\Sym_n}$ is 
a polynomial algebra $\mathbb F[s_1, \dots, s_n]$, where 
$s_1 = P_1$, $s_2 = P_{1, 1}$, etc., are the elementary
symmetric polynomials in $x_1, \dots, x_n$. 
$I_{n, 0}$ is clearly the maximal ideal of $R_n$ 
generated by the regular sequence $s_1, \dots, s_n$. 
Thus Theorem~\ref{thm2}(b) fails if $d = 0$.

Now suppose $d = 1$ and $n = 2p - 1$, where
$\Char(\mathbb F) = p$. By Theorem~\ref{thm1}(b), $I_{n, 1}$
is generated by the $n + 1$ elements $P_2, \dots, P_{n-1}, P_{n+1}$ and
$P_{\underbrace{\text{\tiny $2, \ldots, 2$}}_{\text{$p$ times}}}$.

Since we are in characteristic $p$, $P_{n+1} = P_{2p} = P_2^p$,
is a redundant generator.  In other words, $I_{n, 1}$
is generated by the $n$ elements $P_2, \dots, P_{n-1}, P_n$ and
$P_{\underbrace{\text{\tiny $2, \ldots, 2$}}_{\text{$p$ times}}}$.
By Lemma~\ref{lem.thm2}(c) these elements form a regular sequence
in $R_n$.  This shows that Theorem~\ref{thm2}(b) fails
for $d = 1$ and $n = 2p-1$. 
\qed
\end{remark}

% \begin{remark}
% After writing this paper it was brought to our attention that
% \ref{thm1}(a) and
% \ref{thm2}(a) are also proved in \cite{CKW}.
% \end{remark}
                   
\section*{Acknowledgments} 
The authors are grateful to L. Avramov, W. Dwyer, G. Lyubeznik, 
B. Sturmfels and J. Weyman for their helpful comments.

\end{document}